\documentclass[a4paper,11pt]{article}

\usepackage{amsmath}
\usepackage{amssymb}
\usepackage{amsfonts}
\usepackage{amsthm}
\usepackage{imakeidx}
\usepackage{hyperref}
\usepackage{enumitem}
\usepackage{mathtools}
\usepackage{tikz-cd}
\usepackage{diagbox}

\newtheorem{definition}{Definition}[section]
\newtheorem{lemma}{Lemma}
\newtheorem{theorem}{Theorem}

\newcommand{\Eins}{{\mathbf 1}}
\newcommand{\CC}{\mathbb C}
\newcommand{\HH}{\mathbb H}
\newcommand{\QQ}{\mathbb Q}
\newcommand{\RR}{\mathbb R}
\newcommand{\ZZ}{\mathbb Z}
\newcommand{\ca}{\mathcal}
\newcommand{\ti}{\tilde}
\newcommand{\ol}{\overline}
\newcommand{\GL}{\mathop{\mathrm{GL}}}
\newcommand{\Sp}{\mathop{\mathrm{Sp}}}
\newcommand{\Gr}{\mathop{\mathrm{Gr}}}
\newcommand{\SL}{\mathop{\mathrm{SL}}}
\newcommand{\iso}{\mathop{\mathrm{iso}}}
\newcommand{\Temb}{\mathop{\mathrm{Temb}}}
\newcommand{\RST}{\mathop{\mathrm{RST}}}

\begin{document}

\title{Kuga Varieties of Polarised Abelian Surfaces}
 \author{Wing Kei Flora Poon}
 \maketitle

\section{Introduction} \label{intro}
An $n$-fold Kuga variety, which we will refer to as a Kuga variety, is a variety over a Siegel modular variety such that each fibre is a product of $n$ copies of the abelian variety or Kummer variety to which it corresponds in the base. In this paper, we study the Kodaira dimension of $n$-fold Kuga varieties, $\ca{X}^n_p$, over the moduli spaces $\ca{A}_p$ of
$(1,p)$-polarised abelian surfaces with canonical level structure for prime $p
\geq 3$. 

There are results about the Kodaira dimensions of some special kinds of Kuga varieties: \cite{v}, \cite{fv} have proved the unirationality of Kuga families, i.e. the $1$-fold Kuga varieties, over moduli spaces of principally polarised abelian varieties of dimension $4$ and $5$; in \cite{psms}, we computed the Kodaira dimension of any Kuga variety over moduli spaces of principally polarised abelian varieties of dimension $g \geq 2$. It is therefore natural to ask for the Kodaira dimension of Kuga varieties of other kinds, for example, the $\ca{X}^n_p$ described above. 

A connection between modular forms and differential forms on arbitrary Kuga varieties is established in \cite{m}. This gives us some information about the Kodaira dimensions, assuming a specific compactification for the Kuga varieties, which is referred to as a Namikawa
compactification in \cite{psms}.
Specifically, \cite[Theorem~1.3]{m} translates to: 

\begin{theorem} \label{Ma's thm}
Let $X$ be a Namikawa compactification of $\ca{X}_p^n$. Then
\[
\kappa(\ol{\ca{A}_p}, (n+3)\ca{L}-\Delta_{\ca{A}}) \leq \kappa(K_X) \leq 3
\] 
where $\ol{\ca{A}_p}$ is a toroidal compactification of $\ca{A}_p$, $\ca{L}$ is the $\QQ$-line bundle of weight $1$ modular forms of
$\Gamma_p$ and $\Delta_{\ca{A}}$ is the boundary divisor of
$\ol{\ca{A}_p}$. 
\end{theorem}

The Kodaira dimension is a birational invariant by definition, so
$\kappa(\ca{X}^n_p) = \kappa(X)$. Moreover, if $X$ has canonical singularities,
then $\kappa(X) = \kappa(K_X)$, and hence
$\kappa((n+3)\ca{L}-\Delta_{\ca{A}})$ is a lower bound for
$\kappa(\ca{X}^n_p)$. 

We say a Kuga variety
$\ca{X}^n_p$ is of \emph{relative general type} if its Kodaira
dimension equals the dimension of the base $\ca{A}_p$ of $\ca{X}^n_p$,
namely $3$, which is also the maximum value $\kappa(\ca{X}^n_p)$ can attain.

This paper is divided into two parts: in Section~\ref{canonical} we show that  for $n > 2$ and any $p$,
the particular Namikawa compactification $X$ of $\ca{X}^n_p$ constructed in
\cite{psms} has canonical singularities; in Section~\ref{cusp form} we search for a lower bound of $(p,n)$ for which
$\kappa((n+3)\ca{L}-\Delta_{\ca{A}}) = 3$. We summarise our result in the
following theorem:
\begin{theorem} \label{main result}
A Kuga variety $\ca{X}^n_p$ is of relative general type if
\begin{itemize}
\item $p \geq 3$ and $n \geq 4$;
\item $p \geq 5$ and $n \geq 3$.
\end{itemize}
\end{theorem}

Combining the results of \cite{gh} and \cite{hs}, which say
$\ca{X}^0_p = \ca{A}_p$ is of general type for $p \geq 37$, we can
mark on the $(p,n)$-plane a region for which the Kuga varieties are of relative general type as in
Figure~1.
\begin{center}
\begin{tikzpicture}
\draw[->] (0,0) -- (0,6);
\draw (0,0) -- (8,0);
\draw [dashed] (8.0,0) -- (9.6,0);
\draw[->] (9.6,0) -- (11.74,0);
\coordinate (Figure) at (6,-1);
\filldraw (Figure) node[] {Figure 1};
\foreach \x in {0,0.8,1.6,2.4,3.2,4,4.8,5.6,6.4,7.2,10.4,11.2}
\draw (\x,2pt) -- (\x,-2pt);
\foreach \y in {0,1,2,3,4,5}
\draw (2pt,\y) -- (-2pt,\y);
\coordinate (n) at (0,6);
\filldraw (n)  node[above] {$n$};
\coordinate (p) at (11.74,0);
\filldraw (p)  node[right] {$p$};
\coordinate (label) at (7,4);
\filldraw (label) [red] node[] {of relative general type};

\coordinate (0) at (0,0);
\filldraw (0)  node[left] {0};
\coordinate (1) at (0,1);
\filldraw (1)  node[left] {1};
\coordinate (2) at (0,2);
\filldraw (2)  node[left] {2};
\coordinate (3) at (0,3);
\filldraw (3)  node[left] {3};
\coordinate (4) at (0,4);
\filldraw (4)  node[left] {4};
\coordinate (5) at (0,5);
\filldraw (5)  node[left] {5};

\coordinate (1) at (0.8,0);
\filldraw (1)  node[below] {1};
\coordinate (2) at (1.6,0);
\filldraw (2)  node[below] {2};
\coordinate (3) at (2.4,0);
\filldraw (3)  node[below] {3};
\coordinate (4) at (3.2,0);

\coordinate (5) at (4,0);
\filldraw (5)  node[below] {5};
\coordinate (6) at (4.8,0);

\coordinate (7) at (5.6,0);
\filldraw (7)  node[below] {7};
\coordinate (8) at (6.4,0);

\coordinate (9) at (7.2,0);

\coordinate (40) at (10.4,0);

\coordinate (41) at (11.2,0);
\filldraw (41)  node[below] {37};

\draw[red]  (4,3) -- (5.6,3);
\draw[red]  (4,4) -- (4,3);
\draw[red]  (2.4,4) -- (4,4);
\draw[red]  (2.4,6) -- (2.4,4);
\draw[red,dashed]  (4.8,3) -- (11.2,3);
\draw[red,dashed]  (11.2,1) -- (11.2,3);
\draw[red]  (11.2,1) -- (11.2,0);
\end{tikzpicture}
\end{center}

\subsection{Construction of Kuga varieties} \label{notations}

The moduli space $\ca{A}_p$ of $(1,p)$-polarised abelian surfaces with canonical level structure, defined and studied in \cite[Chapter I.1]{hkw1}, is given as the quotient of the Siegel upper half plane $\HH_2$ by the action of a certain arithmetic subgroup $\Gamma_p$ of $\Sp(4,\ZZ)$.
Recall the Siegel upper half space of degree $g$ is defined as 
\[
\HH_g = \left \{\tau
\in M_{g \times g}(\CC): \tau = \prescript{t}{}{\tau}, \mathop{\mathrm{Im}}\tau > 0 \right \}.
\] 
Note that for any integer $k \geq 2$, the Grassmannian
$\Gr (2,\CC^k)$
is isomorphic to the orbit space $M_{k\times 2}(\CC)/\GL(2,\CC)$
of all $k \times 2$ matrices modulo right multiplication by the invertible matrices in
 $\GL(2,\CC)$. So the Siegel upper half plane $\HH_2$ can be identified with a subset
of $\Gr (2, \CC^4)$ by sending an element $\tau$ to the
$\GL(2,\CC)$-equivalence class of block matrices:
\[
\tau \mapsto 
\begin{bmatrix}
\tau\\
\Eins_2
\end{bmatrix}.
\]

For any prime $p \geq 3$, we define 
\[
\Gamma_p = \left\{ \gamma \in \Sp(4,\ZZ): \gamma - \Eins \in 
\begin{pmatrix}
\ZZ & \ZZ & \ZZ & p\ZZ\\
p\ZZ & p\ZZ & p\ZZ & p^2\ZZ\\
\ZZ & \ZZ & \ZZ & p\ZZ\\
\ZZ & \ZZ & \ZZ & p\ZZ\\
\end{pmatrix}\right\},
\]
which is of finite index in $\Sp(4,\ZZ)$.
This group $\Gamma_p$ acts on $\HH_2 \subset \Gr(2,\CC^4)$ by left multiplication on the GL$(2,\CC)$-equivalence classes of block
matrices, which is an analogue of linear fractional transformations:
\[
\gamma \cdot \tau = 
\begin{bmatrix}
\gamma \cdot 
\begin{pmatrix}
\tau \\ \Eins_2
\end{pmatrix}
\end{bmatrix},
\quad\gamma \in \Gamma_p, \tau \in \HH_g.
\]
The quotient of $\HH_2$ by this action of $\Gamma_p$ gives the moduli space
$\ca{A}_p$.

In \cite[Chapter I.2B]{hkw1}, a $1$-fold Kuga variety $\ca{X}^1_p$ over $\ca{A}_p$ is constructed. This method can be extended to construct an $n$-fold Kuga variety $\ca{X}^n_p$ over $\ca{A}_p$ for any positive integer $n$ by  defining an extension $\ti{\Gamma}^n_p$ of $\Gamma_p$ and a left action of it on $\CC^{2n}\times \HH_2$ which descends to that of $\Gamma_p$ on $\HH_2$. First, by identifying
$\CC^{2n}$ with the set of $n\times 2$ complex matrices, we can
identify $\CC^{2n} \times \HH_2$ with a subset of $\Gr
(2,\CC^{n+4})$ by sending an element $(Z,\tau)$ to a $\GL(2,\CC)$-equivalence class of block matrices:
\[
(Z,\tau)\mapsto 
\begin{bmatrix}
Z\\
\tau\\
\Eins_2\\
\end{bmatrix}.
\]

We define the following group 
\[
\ti{\Gamma}_p^n = 
\left\{ (l,\gamma) =
\begin{pmatrix}
  \Eins_{n} & l \\
  0 & \gamma
\end{pmatrix}
\in \text{M}_{n \times 4}(\ZZ)\rtimes \Gamma_p
: 
\begin{matrix}
\gamma \in \Gamma_p
\end{matrix}
 \right\}.
\]

The group $\ti{\Gamma}_p^n$ acts on $\CC^n \times \HH_2$ by
left multiplication on the $\GL(2,\CC)$-equivalence classes of block
matrices. Explicitly, for $\ti{\gamma} = (l,\gamma) \in
\ti{\Gamma}^n_p$ and $\ti{\tau} = (Z,\tau) \in \CC^{2n} \times
\HH_2$, then
\begin{equation}
\ti{\gamma} \cdot \ti{\tau} = 
\begin{bmatrix}
\begin{pmatrix}
  \Eins_{n} & l \\
  0 & \gamma
\end{pmatrix} \cdot 
\begin{pmatrix}
Z \\ \tau \\ \Eins_2
\end{pmatrix}
\end{bmatrix} = 
\begin{bmatrix}
Z + l \cdot \tau\\
\gamma \cdot \begin{pmatrix}\tau \\ \Eins_2\\ \end{pmatrix}\\
\end{bmatrix} = 
\begin{bmatrix}
(Z + l \cdot \tau)\cdot N\\
\gamma \cdot \tau\\
\Eins_2\\
\end{bmatrix}
\label{eq:1} \tag{$*$}
\end{equation}
for some $N \in \GL(2,\CC)$. 

The quotient of $(\CC^{2n}
\times \HH_2)$ by this action of $\ti{\Gamma}^n_p$ gives the Kuga variety $\ca{X}_p^n$. The projection $\CC^{2n} \times \HH_2
\rightarrow \HH_2$ induces a map $\pi \colon \ca{X}^n_p \rightarrow \ca{A}_p$. Indeed \cite{psms} \cite{n}, all fibres are the product of $n$ copies of the torus parametrised by the corresponding point in the base $\ca{A}_p$, up to a change of basis of $\CC^2$ for each copy.

\section{Canonical singularities}\label{canonical}

\subsection{Namikawa compactification} \label{Namikawa}

The explicit construction of a Namikawa compactification $X$ of the Kuga variety $\ca{X}^n_p$ is the main purpose of section 1 in \cite{psms}. Here we briefly introduce its idea and some notations to be used in the latter sections.

\begin{definition} \label{Nam def}
A \emph{Namikawa compactification} of $\ca{X}_p^n$ is an irreducible
normal projective variety $X$ containing $\ca{X}_p^n$ as an open
subset, together with a projective toroidal compactification
$\ol{\ca{A}_p}$ of $\ca{A}_p$ for which the following conditions hold.
\begin{enumerate}
  \item $\pi \colon \ca{X}_p^n\to \ca{A}_p$ extends to a projective morphism
    $\ol{\pi}\colon X\to \ol{\ca{A}_p}$;
  \item every irreducible component of
    $\Delta_X:=X\smallsetminus \ca{X}_p^n$ dominates an
    irreducible component of
    $\Delta_\ca{A}:=\ol{\ca{A}_p}\smallsetminus\ca{A}_p$.
\end{enumerate}
\end{definition}

Therefore $X$ sits inside the commutative diagram
\[
\begin{tikzcd}
\ca{X}^n_p \arrow[hookrightarrow]{r}{}\arrow{d}{\pi} 
  & X  \arrow{d}{\ol{\pi}} \\
\ca{A}_p \arrow[hookrightarrow]{r}[swap]{} 
  & \ol{\ca{A}_p}
\end{tikzcd}
\]
and $\ol{\pi}$ does not contract any divisor.

Namikawa compactifications are constructed by toroidal methods in \cite{n}. We now describe the partial compactification at an integral boundary component $\ti{F}$ of rank $g' \leq 2$ of $\ca{X}^n_p$ which leads to a Namikawa compactification. The boundary component $\ti{F}$ can be written as $\CC^n \times F$, an extension of a rank $g'$ boundary component $F < \HH_2$ of $\ca{A}_p$, where $F$ corresponds to a rank $g'' := g-g'$ isotropic sublattice of $\ZZ^4$ defined up to the transitive action of $\Sp(4,\ZZ)$. Let $\ti{\ca{P}}(\ti{F})$ be the stabiliser subgroup of $\ti{F}$ in $\RR^{4n} \rtimes \Sp(4,\RR)$, which can be embedded in $\GL(n+4,\RR)$. We also define the following subgroups of $\ti{\ca{P}}(\ti{F})$:
\begin{align*}
    \ti{\ca{P}}'(\ti{F}) &:= \text{ Centre of the unipotent radical of } \ti{\ca{P}}(\ti{F}) \\
    \ti{\Upsilon}^n &:= \ti{\ca{P}}'(\ti{F}) \cap \ti{\Gamma}^n_p\\
    \ti{P''}(\ti{F}) &:= \big(\ti{\ca{P}}(\ti{F}) \cap \ti{\Gamma}^n_p\big)/\ti{\Upsilon}^n
\end{align*}
An explicit expression of the above matrix groups can be found in \cite[Section~2]{n} and \cite[Section~1]{psms}. Also note that both $\ti{\Upsilon}^n$ and $\ti{P''}(\ti{F})$ inherit the
action of $\ti{\Gamma}^n_p$ on $\CC^{2n} \times \HH_2$.

Consider the Siegel domain realisation of $\HH_2$, which gives $\CC^{2n} \times \HH_2$ as an open subset of
\[
(\CC^{g'n} \times \CC^{g''n}) \times (\HH_{g'} \times
M_{g'\times g''}(\CC) \times M_{g''\times g''}^{sym}(\CC)).
\]
Then taking the partial quotient by $\ti{\Upsilon}^n$ near the boundary
component $\ti{F}$ corresponds to translations in the
imaginary directions of the factors $\CC^{g''n}$ and
$M_{g''\times g''}^{sym}(\CC)$ respectively. Therefore there is an open subset of
\[
\CC^{g'n} \times (\CC^*)^{g''n} \times \HH_{g'} \times
\CC^{g'\times g''} \times (\CC^*)^{g''\times g''}_{sym}
\]
that uniformises $\ca{X}^n_p$.

Using a suitable cone decomposition $\Sigma (\ti{F})$, we can extend
the action of $\ti{P''}(\ti{F})$ to a smooth torus embedding
$\Temb(\Sigma (\ti{F}))$ for the torus part
$(\CC^*)^{g''n} \times (\CC^*)^{g''\times g''}_{sym}$. The
partial compactification of $\ca{X}^n_p$ at the boundary component
$\ti{F}$ is then given as the quotient of (an open subset of) the torus bundle
\[
\ti{X}(\ti{F}):= \HH_{g'} \times \CC^{g'g''} \times \CC^{g'n}
\times \Temb(\Sigma (\ti{F}))
\]
by the action of $\ti{P''}(\ti{F})$. Note that this decomposition of $\ti{X}(\ti{F})$ into its factors is preserved by the quotient.

In practice, such a cone decomposition $\Sigma(\ti{F})$ can be given by an extension of the perfect cone decomposition near the boundary component $F$ of $\ca{A}_p$. Furthermore, this extension can be chosen carefully to satisfy more conditions as listed in \cite[proposition~1.4]{psms}. In particular, the set of cone decompositions is compatible at each cusp such that it results in a Namikawa compactification $X$. Also, the local uniformising space $\ti{X}(\ti{F})$ of $X$ has canonical singularities.

We will prove in the remaining subsections that for any $p$ and $n
> 2$, this Namikawa compactification $X$ of $\ca{X}^n_p$ has
canonical singularities.

\subsection{The general strategy} \label{strategy}

We will separately examine the singularities in the interior and the boundary of $X$, and check if they are canonical by applying the Reid--Shepherd-Barron--Tai (RST) criterion.

We will need the following set up to state the RST criterion \cite{r}: 
Suppose $G$ is a finite group acting on the complex vector space $\CC^m$ linearly. For a non-trivial element
$\gamma \in G$ of order $k$, the eigenvalues of
its action on $\CC^m$ can be expressed as an $m$-tuple $(\xi^{\alpha_1},
\cdots, \xi^{\alpha_m})$, with $\xi$ being a primitive $k$-th root of unity and $\alpha_j$ being a non-negative integer less than $k$ for any $j$. We define, with dependence on the choice of  $\xi$, the \emph{type} of $\gamma$ to be 
\[\frac{1}{k}(\alpha_1,\cdots, \alpha_m)\]
and its associated \emph{RST sum} to be
\[
\RST(\gamma) := \sum_{i = 1}^{m} \frac{\alpha_i}{k}.
\]
Furthermore, we say that $\gamma$ is a \emph{quasi-reflection} if all but one $\alpha_j$'s are
$0$, or equivalently $\gamma$ preserves a divisor. 

The RST criterion is then given by the following:
\begin{theorem}[{\cite[4.11]{r}}]
    Let $G$ be a finite group which acts on $\CC^m$ as above. Then $\CC^m/G$ has a canonical singularity if $G$ contains no quasi-reflection, and if every non-trivial element $\gamma \in G$ satisfies the inequality 
    \[\RST(\gamma) \geq 1.\]
\end{theorem}

Note that, since we need to check the above inequality involving the RST sum for every element in $G$, it does not matter which root of unity $\xi$ was chosen to give the type of a generator $\gamma$ of $G$.

\subsubsection{Strategy in the interior}\label{strategyInt}
In the interior $\ca{X}^n_p$ of $X$, a singularity corresponds to a point $\ti{\tau} = (Z, \tau)$ in $\CC^{2n} \times \HH_2$ fixed by $\ti{\Gamma}^n_p$. We are allowed to apply the RST criterion to check if $\ti{\tau}$ corresponds to a canonical singularity:
suppose $\ti{\gamma}$ is an element in the isotropy group $\iso(\ti{\tau}) < \ti{\Gamma}^n_p$ of $\ti{\tau}$. By \eqref{eq:1}, it is clear that one can consider the action of $\ti{\gamma}$ separately as the action of $\gamma$ on the $\HH^2$ factor and that of $\ti{\gamma}$ on the $\CC^{2n}$ factor. Also, $\ti{\gamma}$ fixes $\ti{\tau}$ only if $\gamma$ fixes $\tau$. 
The isotropy group $\iso(\ti{\tau})$ of $\ti{\tau}$ in $\ti{\Gamma}^n_p$ is finite, so any nontrivial element $\ti{\gamma} = (l,\gamma)$ in $\iso(\ti{\tau})$ is a torsion element and $l = 0$. As a result of \cite[Theorem~4.1]{t}, the
induced action of any element $\gamma \in \iso(\tau) \leq \Gamma_p$ of
order $k$ on the tangent space $T_\tau(\HH^2)$ can be diagonalised
under suitable local coordinates. It will be shown that $\ti{\gamma}$
also acts diagonally on $T_Z(\CC^{2n})$. This gives us the finite dimensional
representation of $\iso(\ti{\tau})$ required for the application of the RST criterion.

Note that it suffices to apply the RST criterion at a limited number of singularities in $\ca{X}^n_p$: 
\begin{lemma}\label{lem1}
    Let $\ti{\tau} = (Z, \tau)$ be a point in $\CC^{2n} \times \HH_2$ that corresponds to a canonical singularity in $\ca{X}^n_p$. Then either $\tau$ corresponds to a canonical singularity in $\ca{A}_p$, or $\iso(\ti{\tau})= \langle \ti{\sigma} := (0, -\Eins_4)\rangle < \ti{\Gamma}^n_p$. In the latter case, $\tau$ corresponds to a smooth point. 
\end{lemma}
\begin{proof}
The isotropy group of $\ti{\tau}$, $\iso(\ti{\tau})$, cannot contain a quasi reflection: according to
\cite[Lemma~7.1]{m}, a non-trivial element $\ti{\gamma} \in \iso(\ti{\tau})$ does not fix any divisor in $\ca{X}^n_p$.

Consider any nontrivial $\ti{\gamma} := (0, \gamma) \in \iso(\ti{\tau})$. If $\gamma$ acts trivially on $\HH_2$, then $\gamma = -\Eins_4$. 

Moreover, by the definition of RST sums, we have
\[\RST(\ti{\gamma}) \geq \RST(\gamma)\]
So $\ti{\tau}$ corresponds to a canonical
singularity in $\ca{X}^n_p$ if $\tau$ corresponds to a canonical singularity in $\ca{A}_p$.
\end{proof}

\subsubsection{Strategy in the boundary}\label{strategyBd}
A singularity in the boundary of $X$ correspond to a point $\ti{\tau}$ in $\ti{X}(\ti{F})$ fixed by $\ti{P''}(\ti{F})$ near a boundary component $\ti{F}$ of rank $g'$. Again, the RST criterion can be applied to check if $\ti{\tau}$ corresponds to a canonical singularity: Let $\ti{\tau}:= (Z, \tau)$, where $Z \in \CC^{2n}$ and $\tau \in \HH_{g'} \times \Temb(\Sigma (\ti{F}))$. As mentioned in section \ref{Namikawa},
$\ti{P''}(\ti{F})$ preserves the decomposition of $\ti{X}(\ti{F})$, so $\ti{\gamma}$ acts on each factors of $\ti{X}(\ti{F})$ separately.
A calculation similar to \eqref{eq:1} shows that locally at $\ti{\tau}$, $\ti{\gamma} = (l,\gamma) \in
\ti{P''}(\ti{F})$ fixes $\ti{\tau}$ only if $\gamma$ fixes
$\tau$. However, different from what we had in section \ref{strategyInt}, $\ti{\gamma}$ may not be a torsion
element, i.e. $l$ could be non-zero. 
Nevertheless, the action of $\ti{\gamma}$ on the tangent space of a resolution of $\ti{X}(\ti{F})$ at $\ti{\tau}$ at $\ti{\tau}$ is of finite order, so
the RST criterion can
be applied there.

The following observations are useful for checking whether these singularities are canonical: 
\begin{enumerate}
        \item \cite[Lemma~1.3]{psms}: Let $(\ti{X}(\ti{F}))^*$ be a smooth $\ti{P''}(\ti{F})$-equivariant resolution of $\ti{X}(\ti{F})$.  If $\ti{P''}(\ti{F})$ has no quasireflection, then the partial compactification $\ti{P''}(\ti{F}) \backslash \ti{X}(\ti{F})$ has canonical singularities if $\ti{P''}(\ti{F}) \backslash (\ti{X}(\ti{F}))^*$ has canonical singularities. In particular, this implies that we can apply the RST criterion at the singularities in $\ti{P''}(\ti{F}) \backslash (\ti{X}(\ti{F}))^*$ instead. 

        \item Let $\ti{\tau} = (Z, \tau)$ correspond to a canonical singularity near $\ti{F}$. Then either $\tau$ corresponds to a canonical singularity in the boundary of $\ol{\ca{A}_p}$,
        or $\iso(\ti{\tau})= \langle \ti{\sigma} := (l, -\Eins_4)\rangle <  \ti{P''}(\ti{F})$ for some $l \in L$. In the latter case, $\tau$ corresponds to a smooth point. The proof is similar to that in Lemma \ref{lem1}. Again, this implies that we only need to apply the RST criterion at a limited number of singularities.
    \end{enumerate}

\subsection{Singularities in the interior of compactification} \label{interior}
In this section, we will identify the singularities in $\ca{X}^n_p$ and show that for $n > 2$, they are all canonical.

First we identify singularities that project to non-canonical
singularities in $\ca{A}_p$. 
It is given in the proof of \cite[Theorem~1.8]{hkw2} that for any odd prime $p$, the singular points in $\ca{A}_p$
are exactly the points that lie on the two disjoint curves $C_1$ and $C_2$. 
Any point on one of these curves corresponds to a point $\tau$ in $\HH_2$, whose isotropy group in $\Gamma_p$ is generated by a single generator. Its induced action the tangent space of $\HH_2$ at $\tau$ is also given there: one can write any point in the tangent space $T_\tau({\HH_2})$ in the form
\[
\begin{pmatrix}
    \tau_1+x & \tau_2+y\\
    \tau_2+y & \tau_3+z
\end{pmatrix}.
\]
So the tuple $(x,y,z)$ can be considered as the local coordinates for $T_\tau{\ca{A}_p}$, and the respective action of a generator of $\iso(\tau)$ on $T_\tau({\HH_2})$ with these coordinates is given by
\begin{align*}
    (x, y, z) &\mapsto (-x, -iy, z) \text{ along } C_1;\\
    (x, y, z) &\mapsto (\rho^2 x, -\rho y, z) \text{ along } C_2 \text{, where } \rho = e^{2 \pi i / 3}.
\end{align*}
Therefore, the chosen generators are of types $\frac{1}{4}(2,3,0)$ and
$\frac{1}{6}(4,5,0)$ when the root of unity $\xi$ is chosen to be $i$ and $e^{2\pi i
  /6}$ respectively on each curve $C_1$ and $C_2$. By applying the RST criterion to the isotropy groups, it is clear that  singularities on $C_1$ are
canonical but those on $C_2$ are not. 

Let $\ti{\tau} := (Z,\tau) \in \CC^{2n} \times \HH_2$ such that $\tau$ corresponds to a point in $C_2$.
Let $\ti{\sigma} := (0, -\Eins_4)$ and $\ti{\gamma} := (0, \gamma)$, where $\gamma$ is the generator of $\iso(\tau)$ with the action on $T_\tau({\HH_2})$ described above. Then either $\iso(\ti{\tau}) = \langle \ti{\gamma} , \ti{\sigma} \rangle$ or $\iso(\ti{\tau}) = \langle \ti{\gamma} \rangle$.

We shall first compute the type of $\ti{\gamma}$. We only need to understand 
the action of $\ti{\gamma}$ at a
point $\ti{Y} = (Z+Y, \tau)$ on the tangent space
$T_{\ti{\tau}}(\CC^{2n} \times \{\tau\}) \simeq T_Z(\CC^{2n})$ to complete the type of $\ti{\gamma}$. To do this, we need the explicit expressions of the set $C_2$ and its isotropy group $\iso(\tau)$ from \cite[Definition~1.5]{hkw2}:
\begin{align*}
        C_2 &=
\left\{
\begin{pmatrix}
  \rho & 0 \\
  0 & \tau_3
\end{pmatrix}
:
\rho = e^{2\pi i /3}, 
\tau_3 \in \HH_1
\right\},
\\
\iso(\tau) &= \left\langle \gamma =
\begin{pmatrix}
 0 & 0 & -1 & 0 \\
 0 & 1 & 0 & 0 \\
 1 & 0 & 1 & 0 \\
 0 & 0 & 0 & 1
\end{pmatrix} \right\rangle.
\end{align*}

Following \eqref{eq:1}, the action of $\ti{\gamma}$ at  $\ti{Y}$ in $T_{\ti{\tau}}(\CC^{2n} \times \{\tau\})$
is given by
\[
\ti{\gamma} \cdot \ti{Y} = 
\begin{bmatrix}
(Z + Y) \cdot N\\
\tau\\
\Eins_2
\end{bmatrix}
\text{, where }
N = 
\begin{pmatrix}
(\rho+1)^{-1} & 0\\
0 & 1
\end{pmatrix}.
\]
Since $\ti{\gamma}$ fixes $(Z, \tau)$, $Z\cdot N = Z$ and
$\ti{\gamma}$ acts on $T_Z(\CC^{2n})$ diagonally by sending the set of local coordinates $Y$ to $Y \cdot N$.

Note that $(\rho + 1)^{-1} = e^{2\pi i \cdot (5/6)}$. So by choosing
the primitive root of unity to be $e^{2\pi i/6}$, which is the same as
that for the $\HH_2$ factor, we have an extra $n$ copies of $5/6$'s and $n$ copies of $0$'s in the RST sum of $\ti{\gamma}$. In other words, the type of $\ti{\gamma}$ is $\frac{1}{6}(4, 5, 0, 5, \cdots , 5, 0, \cdots, 0)$.

As for the type of $\ti{\sigma}$, since $\ti{\sigma}$ acts trivially on $T_{\tau}(\{Z\} \times \HH_2)$, 
the first entries in the type of $\ti{\sigma}$ which correspond to the $\HH_2$ factor are all $0$'s. On the other hand, the calculation in \eqref{eq:1} shows that $\ti{\sigma}$ acts on
the set of local coordinates in $T_Z(\CC^{2n})$ diagonally by $X
\mapsto - X$. So the type of $\ti{\sigma}$ is
$\frac{1}{2}(0, 0, 0, 1, ..., 1)$ when the primitive root of unity
$\xi$ is chosen to be $ -1$. 

Since $\ti{\sigma}$ commutes with $\ti{\gamma}$, we
can draw the following table which shows the type of a non-trivial element $\ti{\gamma}^{k_1}\ti{\sigma}^{k_2} \in \iso(\ti{\tau})$, where $0 \leq k_1 \leq 5$ and $0 \leq k_2 \leq 1$.
\begin{center}
\setlength{\extrarowheight}{.3em}
\begin{tabular}{ c|c|c } 
 \diagbox[width=3em]{$k_1$}{$k_2$} & 0 & 1\\ 
 \hline
 0 & N/A & $\frac{1}{2}(0, 0, 0, 1, \cdots , 1, 1, \cdots , 1)$\\
 1 & $\frac{1}{6}(4, 5, 0, 5, \cdots , 5, 0, \cdots, 0)$ & $\frac{1}{6}(4, 5, 0, 2, \cdots , 2, 3, \cdots , 3)$\\
 2 & $\frac{1}{6}(2, 4, 0, 4, \cdots , 4, 0, \cdots, 0)$ & $\frac{1}{6}(2, 4, 0, 1, \cdots , 1, 3, \cdots , 3)$\\ 
 3 & $\frac{1}{6}(0, 3, 0, 3, \cdots , 3, 0, \cdots, 0)$ & $\frac{1}{6}(0, 3, 0, 0, \cdots , 0, 3, \cdots , 3)$ \\ 
 4 & $\frac{1}{6}(4, 2, 0, 2, \cdots , 2, 0, \cdots, 0)$ & $\frac{1}{6}(4, 2, 0, 5, \cdots , 5, 3, \cdots , 3)$\\
 5 & $\frac{1}{6}(2, 1, 0, 1, \cdots , 1, 0, \cdots, 0)$ & $\frac{1}{6}(2, 1, 0, 4, \cdots , 4, 3, \cdots , 3)$ \\ 
\end{tabular}
\end{center}
The types of all non-trivial elements in $\langle \ti{\gamma}\rangle$ are given by the first column of the table, while that in $\langle \ti{\gamma}, \ti{\sigma} \rangle$ are given by the entire table. 
 Notice the RST criterion only fails when $n \leq 2$:
 \[\RST(\ti{\gamma}^5) < 1.\]
We conclude that for $n > 2$, both $\langle \ti{\gamma}\rangle$ and $\langle \ti{\gamma}, \ti{\sigma} \rangle$ satisfy the RST criterion, and therefore $\ti{\tau}$ is a canonical singularity in $\ca{X}^n_p$, no matter $\iso(\ti{\tau}) = \langle \ti{\gamma} , \ti{\sigma} \rangle$ or $\iso(\ti{\tau}) = \langle \ti{\gamma} \rangle$.

Finally, for any singularity that corresponds to a point in $\CC^{2n} \times \HH_2$ whose isotropy group is  $\langle \ti{\sigma}\rangle$, we only need study the first row of the table: there is no quasi-reflection and the RST inequality is satisfied for any $n$. Therefore such singularity is always canonical.

\subsection{Singularities in the boundary of compactification} \label{boundary}
In this section we will check that every singularity in the boundary of $X $ is canonical. 

First, we identify all the non-canonical singularities in $\ol{\ca{A}_p}$. Consider the compact curves $C_1^*$ and $C_2^*$ containing $C_1$ and $C_2$ in $\ol{\ca{A}_p}$. Then from \cite[Propositions~2.15 and 3.4]{hkw2}, for any odd prime $p$, the complement $\ol{\ca{A}_p} \smallsetminus (C_1^* \cup C_2^*)$ contains only isolated singularities. The types of a generator in the respective isotropy groups are given as $\frac{1}{2}(1, 1, 1)$ or $\frac{1}{3}(1, 2, 1)$. So both isotropy groups satisfy the RST criterion, and these singularities in $X$ are canonical. Therefore, any non-canonical singularity in $X$ has to project down to $C_1^* \smallsetminus C_1$ or $C_2^* \smallsetminus C_2$. 

From the same source above, each set $C_1^* \smallsetminus C_1$ and $C_2^* \smallsetminus C_2$ consists of $(p^2-1)/2$ points, one in each of the rank $1$ boundary components called peripheral components \cite[Definition I.3.105]{hkw1}. \cite[Proposition~2.8]{hkw2} further says that near one of these boundary component
$\ti{F}$, the singularities in $C_1^*$ and $C_2^*$ are represented
by $Q_1 = (i, 0, 0)$ and $Q_2 = (\rho, 0, 0)$ as points in $\HH_1 \times \CC \times \CC$, the Siegel domain realisation of $\HH_2$, with $\rho = e^{2\pi i/3}$.

First consider the singularity in $X$ associated to $Q_2$: let $\ti{\tau} := (Z, \tau) \in \ti{X}(\ti{F})$ such that $\tau = Q_2$. From \cite[Propositions~2.5 and 2.8]{hkw2}, the stabiliser
subgroup of $\tau$ in $P''(F) \simeq (\ca{P}(F) \cap \Gamma_p) /
(\ca{P'}(F) \cap \Gamma_p)$ is generated by the order $6$
element
\[
\gamma =
\begin{pmatrix}
0 & 0 & -1 & 0\\
0 & 1 & 0 & 0\\
1 & 0 & 1 & 0\\
0 & 0 & 0 & 1
\end{pmatrix}.
\]
Let $\ti{\gamma} := (l, \gamma)$ be the corresponding generator in $\iso(\ti{\tau})$, and let $\ti{\sigma} := (l, -\Eins_4)$ for some $l \in L$. Then again $\iso(\ti{\tau}) = \langle \ti{\gamma}\rangle$ or  $\iso(\ti{\tau}) = \langle \ti{\gamma}, \ti{\sigma} \rangle$.
To find the types of elements in $\iso(\ti{\tau})$, we consider their actions on each factor of $\ti{X}(\ti{F})^*$, a $\ti{P''}(\ti{F})$-equivariant resolution of $\ti{X}(\ti{F})$.
\cite[Lemmas~5.1 and 5.2]{t}
describes such a resolution of singularities for the moduli space of polarised abelian $g$-folds, as well as a formula for the RST sum of a generator $\gamma$ in the isotropy group. Explicitly when $g = 2$, there are three factors in the resolution of $\HH_2$: $\HH_{g'}, \CC^{g'g''}$ and a torus at infinity. The following submatrices are extracted from
the entries $\gamma_{ij}$ of $\gamma$:
\[
\gamma' = 
\begin{pmatrix}
\gamma_{11} & \gamma_{13}\\
\gamma_{31} & \gamma_{33}
\end{pmatrix}, \qquad
U = 
\begin{pmatrix}
\gamma_{22}
\end{pmatrix}.
\]
Suppose $\gamma'$ has eigenvalues $\lambda^{\pm1}$ and $U$ has eigenvalue $\mu$. Then the eigenvalues of the action of $\gamma$ on the tangent space of the $\HH_{g'}$ factor, the $\CC^{g'g''}$ factor and the torus at infinity in the resolution of $\HH_2$ are $\lambda$, $\lambda\mu$ and $0$ respectively.

In our case, $\gamma'$ has eigenvalues $e^{\pm 2\pi i /6}$ and $U$ has eigenvalue $1$. Therefore, when $e^{2\pi i /6}$ is chosen to be the primitive
root of unity, the $\HH_{g'}$, $\CC^{g'g''}$ and the torus at infinity factors contribute $\frac{2}{6}$, $\frac{1}{6}$ and
$0$ to the RST sum respectively.

For the RST sum over the remaining $\CC^n \times (\CC^*)^n$ factor of $\ti{\tau}$, follow \eqref{eq:1} and consider the action of $\ti{\gamma}$ on $\ti{Y} := (Z+Y, \tau)$ in the tangent space
at $Z$ of the resolved $\CC^{2n}$ factor:
\[
\ti{\gamma} \cdot \ti{Y} = 
\begin{bmatrix}
\left (Z'+Y \right ) \cdot N\\
\tau\\
\Eins_2
\end{bmatrix}
\text{ where }
Z' = Z+l\cdot \begin{pmatrix}
\tau \\ \Eins_2
\end{pmatrix} \text{ and }
N = 
\begin{pmatrix}
\frac{1}{\rho + 1} & 0\\
0 & 1
\end{pmatrix}.\]
Again $\ti{\gamma}$ fixes $\ti{\tau}$, so $Z' \cdot N = Z$ and
$\ti{\gamma}$ acts on the tangent space diagonally by sending the local coordinates $Y$ to $Y \cdot N$.  The eigenvalues of the action are the
eigenvalues of $N$, which are $e^{2\pi i \cdot (5/6)}$ and $1$. When
we choose $e^{ 2\pi i /6}$ to be the primitive root of unity for
$\ti{\gamma}$, which is the same choice as the other factors, they contribute $n$ copies of $\frac{5}{6}$ and $n$
copies of $0$ to the RST sum.

Do the same for $\ti{\sigma}$ to find $\RST(\ti{\sigma})$: write $\sigma = -\Eins_4$ and
consider the submatrices $\sigma'$ and $U$ extracted from $\sigma$ in
the same way as above. Their eigenvalues are $\{-1, -1\}$ and $-1$
respectively, which contribute $0$ to the RST sum for all $3$ factors of the solution of $\HH_2$ after resolving. Following \eqref{eq:1}, the action of $\ti{\sigma}$
on the tangent space of the resolved $\CC^n \times (\CC^*)^n$ factor at $Z$ is again multiplication by $-1$ to the local coordinates $Y$. 

Therefore, we can draw a similar
table as in the previous subsection for each element
$\ti{\gamma}^{k_1}\ti{\sigma}^{k_2} \in \iso(\ti{\tau})$, where $0 \leq k_1 \leq 5$ and $0 \leq k_2 \leq 1$:
\begin{center}
\setlength{\extrarowheight}{.3em}
\begin{tabular}{ c|c|c } 
 \diagbox[width=3em]{$k_1$}{$k_2$} & 0 & 1\\ 
 \hline
 0 & N/A & $\frac{1}{2}(0, 0, 0, 1, \cdots , 1, 1, \cdots , 1)$\\
 1 & $\frac{1}{6}(2, 1, 0, 5, \cdots ,5, 0, \cdots, 0)$ & $\frac{1}{6}(2, 1, 0, 2, \cdots , 2, 3, \cdots , 3)$\\
 2 & $\frac{1}{6}(4, 2, 0, 4, \cdots , 4, 0, \cdots, 0)$ & $\frac{1}{6}(4, 2, 0, 1, \cdots , 1, 3, \cdots , 3)$\\ 
 3 & $\frac{1}{6}(0, 3, 0, 3, \cdots , 3, 0, \cdots, 0)$ & $\frac{1}{6}(0, 3, 0, 0, \cdots , 0,3, \cdots , 3)$ \\ 
 4 & $\frac{1}{6}(2, 4, 0, 2, \cdots , 2, 0, \cdots, 0)$ & $\frac{1}{6}(2, 4, 0, 5, \cdots , 5, 3, \cdots , 3)$\\
 5 & $\frac{1}{6}(4, 5, 0, 1, \cdots , 1, 0, \cdots, 0)$ & $\frac{1}{6}(4, 5, 0, 4, \cdots , 4, 3, \cdots , 3)$ \\ 
\end{tabular}
\end{center}
One can check that there is no quasi-reflection, and the RST sum is at least $1$ everywhere on the table. So the RST criterion is satisfied for both $\langle \ti{\gamma}\rangle$ and $\langle \ti{\gamma}, \ti{\sigma}\rangle$. Thus for
all $n \geq 1$, the singularity in $X$ that corresponds to $(Z,Q_2)$ is canonical.

Now we replace $Q_2$ by $Q_1$ everywhere in the above to check whether
the other singularity in the boundary component $\ti{F}$ is canonical or not. Again, let $\ti{\tau} = (Z, \tau)$ such that $\tau = Q_1$. The stabiliser
subgroup of $\tau = Q_1$ in $P''(F)$ is generated by the order $4$ element
\[
\gamma =
\begin{pmatrix}
0 & 0 & -1 & 0\\
0 & 1 & 0 & 0\\
1 & 0 & 0 & 0\\
0 & 0 & 0 & 1
\end{pmatrix}.
\]
First we caluclate $\RST(\gamma)$. Extract the submatrices $\gamma'$ and $U$ as before. The eigenvalues
of $\gamma'$ are $\pm i$ and the eigenalue of $U$ is $1$. When $i$ is
the chosen primitive root of unity, the $\HH_{g'}$ factor, the $\CC^{g'g''}$ factor and the torus at infinity in the resolution of $\HH_2$ contribute a $\frac{2}{4}$, a $\frac{1}{4}$ and a $0$ to the RST sum respectively.
Consider the action of $\ti{\gamma}$ at $\ti{Y} = (Z+Y, \tau)$. Then
$\eqref{eq:1}$ gives:
\[
\ti{\gamma} \cdot \ti{Y} = 
\begin{bmatrix}
\left (Z'+Y \right ) \cdot N\\
\tau\\
\Eins_2
\end{bmatrix}
\text{ where }
Z' = Z+l\cdot \begin{pmatrix}
\tau \\ \Eins_2
\end{pmatrix} \text{ and }
N = 
\begin{pmatrix}
-i & 0\\
0 & 1
\end{pmatrix}.\]
Once more $Z' \cdot N = Z$ and $\ti{\gamma}$ acts on the tangent space
diagonally by sending the local coordinates $Y$ to $Y\cdot N$. This action has eigenvalues
$e^{2\pi i \cdot (3/4)}$ and $1$, which contribute $n$ copies of
$\frac{3}{4}$ and $n$ copies of $0$ to the RST sum over the resolved $\CC^n \times (\CC^*)^n$ factor when the primitive root of unity chosen is
$i$.

The RST sums of $\ti{\sigma}$ restricted to each factor is the same as
the case of $Q_2$. 

Therefore we can draw the table for the type of
$\ti{\gamma}^{k_1}\ti{\sigma}^{k_2}\in \iso(\ti{\tau})$, where $0 \leq k_1 \leq 3$ and $0 \leq k_2 \leq 1$:
\begin{center}
\setlength{\extrarowheight}{.3em}
\begin{tabular}{ c|c|c } 
 \diagbox[width=3em]{$k_1$}{$k_2$} & 0 & 1\\ 
 \hline
 0 & N/A & $\frac{1}{2}(0, 0, 0, 1, \cdots , 1, 1, \cdots , 1)$\\
 1 & $\frac{1}{4}(2, 1, 0, 3, \cdots , 3, 0, \cdots, 0)$ & $\frac{1}{4}(2, 1, 0, 1, \cdots , 1, 2, \cdots , 2)$\\
 2 & $\frac{1}{4}(0, 2, 0, 2, \cdots , 2, 0, \cdots, 0)$ & $\frac{1}{4}(0, 2, 0, 0, \cdots , 0, 2, \cdots , 2)$\\ 
 3 & $\frac{1}{4}(2, 3, 0, 1, \cdots , 1, 0, \cdots, 0)$ & $\frac{1}{4}(2, 3, 0, 3, \cdots , 3, 2, \cdots , 2)$ \\ 
\end{tabular}
\end{center}
The RST criterion is satisfied for both $\langle \ti{\gamma}\rangle$ and $\langle \ti{\gamma}, \ti{\sigma}\rangle$, so for
all $n \geq 1$, the singularity in $X$ that corresponds to $(Z,Q_1)$ is canonical.

We summarise our findings in the
following theorem:

\begin{theorem} \label{sing. thm}
Singularities in the Namikawa compactification $X$ of $\ca{X}^n_p$ are
canonical for $n\geq 3$. For $n = 1, 2$, the set of non-canonical
singularities in $X$ is exactly the preimage under $\ol{\pi}$ of the
curve $C_2$ in $\ca{X}^1_p$ and $\ca{X}^2_p$ respectively.
\end{theorem}

\section{Low weight cusp form trick} \label{cusp form}
In this section, we will prove the following theorem:
\begin{theorem} \label{pn thm}
The equality $\kappa(\ol{\ca{A}_p}, (n+3)L-\Delta_{\ca{A}}) = 3$ is
satisfied for the following values of $n$ and $p$:
\begin{itemize}
\item $p \geq 3$ and $n \geq 4$;
\item $p \geq 5$ and $n \geq 3$.
\end{itemize}
\end{theorem}

To find a lower bound for $\kappa((n+3)\ca{L}-\Delta_{\ca{A}})$, which is
the rate of growth with respect to $m$ of the dimension of the space
of weight $m(n+3)$-cusp forms of $\Gamma_p$, we use the ``low weight
cusp form trick", which has been used in this context in \cite{gh} and \cite{gs}, and more widely thereafter.

Suppose $n>N$ and there exists a non-zero weight $3+N$ cusp form $F$
of $\Gamma_p$, that is, $F \in H^0((3+N)\ca{L}-\Delta_A)$. For any non-zero
$F' \in H^0(m(n-N)\ca{L})$,
\[
F^m F' \in H^0(m((n+3)\ca{L}-\Delta_A)).
\]
Fixing $F$, the space of cusp forms in the form of $F^m F'$ then grows
at the same rate as $H^0(m(n-N)\ca{L})$ with respect to $m$, which is known
to be $\ca{O}(m^3)$. So $\kappa((n+3)\ca{L} - \Delta_A) \geq 3$.

Therefore, $\ca{X}^n_p$ is of relative general type if
\[
h : = \dim H^0((3+N)\ca{L}-\Delta_A) > 0.
\]
  To find a lower bound for $h$,
we apply Gritsenko's lifting of Jacobi cusp forms mentioned in \cite[Theorem~3]{g}, which states the existence of an injective lifting
\[
J^{cusp}(k,p) \hookrightarrow S_{k}(\Gamma[p])
\]
where $J^{cusp}(k,p)$ is the space of Jacobi cusp forms of weight $k$
and index $p \geq 1$, and $S_{k}(\Gamma[p])$ is the space of weight
$k$ cusps forms of $\Gamma[p]$, with the paramodular group $\Gamma[p]$
defining the moduli space of $(1,p)$-polarised abelian surfaces
without level structure as $\Gamma[p] \backslash \HH_2$. But since
$\Gamma_p \leq \Gamma[p]$, the image of the lifting is also contained
in $S_{k}(\Gamma_p)$.

From \cite{ez}, dim $J^{cusp}(k,p) \geq j(k,p)$ (equality holds when
$k \geq p$), where
\[
j(k,p) := 
\begin{cases}
    \sum\limits_{j=0}^t \left (\dim M_{k+2j} - \left ( \left \lfloor \frac{j^2}{4p} \right \rfloor + 1 \right )\right ),& \text{if } k \text{ is even}\\
    \sum\limits_{j=1}^{t-1} \left (\dim M_{k+2j-1} - \left ( \left \lfloor \frac{j^2}{4p} \right \rfloor + 1 \right )\right ),& \text{if } k \text{ is odd}
\end{cases}
\]
with $M_r$ being the space of modular forms of weight $r$ for $\SL(2, \ZZ)$.

It is a general fact that 
\[
\dim M_r = 
\begin{cases}
    \left \lfloor \frac{r}{12} \right \rfloor,& \text{if } r \equiv 2 \mod 12\\
    \left \lfloor \frac{r}{12} \right \rfloor + 1,& \text{otherwise}
\end{cases}
\]

By a simple computation, it can be found that the first prime $p$ such
that $j(k,p) > 0$ for $k = 5$ and $6$ are $p = 5$ and $3$
respectively. Note:
\begin{enumerate}
    \item  $\dim (S_k(\Gamma_p)) \geq j(k,p)$ for any $k, p$;
    \item  $j(k,p)$ increases with $p$; 
    \item the isomorphism in \cite[Theorem~1.1]{m}
      \[
      \bigoplus_{m\geq 0}H^0(\ca{X}^n_p,K_{\ca{X}^n_p}^{\otimes m}) =
      \bigoplus_{m\geq 0}M_{(n+3)m}(\Gamma_p)
      \]
    establishes that $\kappa(\ca{X}^n_p)$ is non-decreasing with
    respect to $n$, because the same is true for
    dim$(M_{(n+3)m}(\Gamma_p))$. 
\end{enumerate}
By letting $k = 3+N = 2+n$, this shows that for the values of $n$ and
$p$ stated in Theorem~\ref{pn thm}, $\dim(S_k(\Gamma_p)) \geq j(k,p) \geq
1$. This concludes our proof for Theorem~\ref{pn thm}.

\section{Possible improvements} \label{future work}
By following \cite{hs} and applying the Riemann-Roch theorem on the
exceptional divisor $E$ of a blow-up at a non-canonical singularity in
$\ca{X}^1_p$, we may be able to improve our boundary at $n = 1$ by
finding two consecutive primes $p'$ and $p''$ with $p' < p''$ such that $\kappa(\ca{X}^1_{p'}) <
\kappa(\ca{X}^1_{p''})$. However, that would involve understanding the
intersection behaviour of divisors on the
$4$-fold $E$, which is expected to be complicated. The low density of
prime numbers near $37$ makes the quest less promising: the estimate for $p'$ we
find by this method may not be smaller than $31$.

There are a few more questions that can be asked: for
example, whether the boundary we have drawn can be improved for $p =
5$ and $p = 3$. The image of Gritsenko's lift is not the entire
$S_k(\Gamma_p)$ or even $S_k(\Gamma[p])$, so we might be able to find a
weight $4$ cusp form with respect to $\Gamma_p$ or $\Gamma[p]$ through
other means which improves the bound at $p=5$, and likewise for
$p=3$. Another question is to calculate $\kappa(\ca{X}^n_p)$ for other
$\ca{X}^n_p$ not of relative general type by considering the slope of
Siegel cusp forms of $\Gamma_p$, which is the ratio
between weight and vanishing order at $\infty$, and to draw a boundary
on the $(p,n)$-plane separating the regions with $\kappa(\ca{X}^n_p) =
-\infty$ and $\kappa(\ca{X}^n_p) \geq 0$. We can also extend the
problem by considering $p = 2$, non-prime $p$, or abelian surfaces
without level structure.

\end{document}